\documentclass[a4paper,12pt,reqno]{amsart}
\usepackage{amsfonts}
\usepackage{amsmath}
\usepackage{amssymb}
\usepackage{lipsum}
\usepackage{mathrsfs}
\usepackage{hyperref}
\usepackage{enumerate}
\numberwithin{equation}{section}

\usepackage{graphicx}
\usepackage{enumitem}

\newtheorem{theorem}{Theorem}[section]
\newtheorem{defi}[theorem]{Definition}

\newtheorem{prop}[theorem]{Proposition}
\newtheorem{lemma}[theorem]{Lemma}
\newtheorem{example}[theorem]{Example}

\newtheorem{quest}[theorem]{Question}

\newenvironment{pff}{\hspace*{-\parindent}{\bf Proof \,}}
{\hfill $\Box$ \vspace*{0.2cm}}

\def\R2n{{\mathbb R}^{2n}}

\def\R2{{\mathbb R}^2}
\def\R2n{{\mathbb R}^{2n}}

\usepackage{fancyhdr}
\begin{document}
	\title[Graphs of continuous functions and fractal dimension]
	{Graphs of continuous functions and fractal dimension}
	
	\author[Manuj Verma]{Manuj Verma}
	\address{
		Manuj Verma:
		\endgraf
		Department of Mathematics
		\endgraf
		Indian Institute of Technology, Delhi, Hauz Khas
		\endgraf
		New Delhi-110016 
		\endgraf
		India
		\endgraf
		{\it E-mail address} {\rm mathmanuj@gmail.com}
	}
	
	\author[Amit Priyadarshi]{Amit Priyadarshi}
	\address{
		Amit Priyadarshi:
		\endgraf
		Department of Mathematics
		\endgraf
		Indian Institute of Technology, Delhi, Hauz Khas
		\endgraf
		New Delhi-110016 
		\endgraf
		India
		\endgraf
		{\it E-mail address} {\rm priyadarshi@maths.iitd.ac.in}
	}

	\thanks{...
	}
	\date{\today}
	
	
	\begin{abstract}
		In this paper, we show that, for any $\beta \in [1,2]$, a given strictly positive real-valued continuous function on $[0,1]$ whose graph has upper box-counting dimension less than or equal to $\beta $ can be decomposed as a product of two real-valued continuous functions on $[0,1]$ whose graphs have upper box-counting dimension equal to $\beta$. We also obtain a formula for the upper box-counting dimension of every element of a ring of polynomials in finite number of continuous functions on $[0,1]$ over the field $\mathbb{R}.$  
		
	\end{abstract}

	\maketitle

	\section{Introduction}  
	Let $X$ be a compact subset of $\mathbb{R}$ and let $C(X)$ denote the space of all real-valued continuous  functions on $X$. Then $C(X)$ is a complete metric space with respect to the sup norm. For any $f\in C(X)$, the graph of $f$ is defined as  
	$$G_{f}(X)= \{(x,f(x)): x\in X \}\subset X \times \mathbb{R}.$$
	For simplicity, we denote the graph of $f$ by $G(f)$ instead of $G_f[0,1]$ when $X=[0,1]$ and $G_f(X)$ for $X\ne [0,1]$.
	Throughout this paper, $\dim_{H}(.)$, $\dim_{P}(.)$ $\dim_{A}(.),\overline\dim_{B}(.)$ and $\underline \dim_{B}(.)$  denote Hausdorff dimension, packing dimension, Assouad dimension, upper box dimension and lower box dimension, respectively. For more details, we refer \cite{Falconer,fraser}.
	\par In 1986, Mauldin and Williams\cite{mauldin1} were the first who tried to decompose a continuous function in terms of fractal dimension. They obtained the following interesting result of decomposition:
	
			For any $f\in C[0,1]$ and $\beta=1$, there exist two continuous functions $g,h\in C[0,1]$ such that  
	\begin{equation}\label{mwr}
	f=g+h~~\text{and}~~ \dim_{H}G(g)=\dim_{H}G(h)= \beta .	
	\end{equation}
	 They have given only the existence of this decomposition. After that in 2000, Wingren \cite{win} gave a techinque for its construction and also introduced the Result \eqref{mwr} in terms of lower box dimension. In 2011, K. J. Falconer and J. M. Fraser \cite{falconer2011horizon} gave a formula for upper box dimension of the graph of sum of two continuous functions, which is as follows: 
	 \begin{theorem}\label{upf} Let $g,h \in C[0,1]$. Then 
		$$\overline{\dim}_BG(g+h)\leq\max \{\overline{\dim}_BG(g), \overline{\dim}_BG(h)\}.$$ 
	\end{theorem}
	 In 2013, F. Bayart and Y. Heurteaux \cite{bayar} verified the Result \eqref{mwr} for $\beta = 2$ and proposed an open problem: Is the Result \eqref{mwr} true for any $\beta\in [1,2] ?$     
	 Later, in 2013, J. Liu and J. Wu \cite{liu01} gave a positive answer to the above question and proved the Result \eqref{mwr} for any $\beta \in  [1,2]$. From this we can say that there is no general formula for Hausdorff dimension of graph of sum of two continuous functions like upper box dimension.
 After that, in 2016, J. Liu, B. Tan and J. Wu \cite{liu02} have given  the packing dimension and upper box dimension version of the Result \eqref{mwr}. Precisely, they gave the following results:     
	\begin{theorem}\cite{liu02}
		Let $\beta \in [1,2]$, and $f \in C[0,1]$. If $\dim_{P}G(f)\leq \beta $, then there are functions $g,h\in C[0,1]$ such that 
		$$f=g+h~~\text{and }~~\dim_{P}G(g)=\dim_{P}G(h)=\beta.$$	
	\end{theorem} 
	\begin{theorem}\label{liu@}\cite{liu02}
		Let $\beta \in [1,2]$, and $f \in C[0,1]$. Then there are functions $g,h\in C[0,1]$ such that 
		$$f=g+h~~\text{and }~~\overline\dim_{B}G(g)=\overline\dim_{B}G(h)=\beta.$$
		if and only if $~~\overline \dim_{B}G(f)\leq \beta.$
	\end{theorem}
	Recently in 2020, J. Liu and D. Liu \cite{liu03} establised  the decomposition Result \eqref{mwr} for lower box dimension and $\beta \in [1,2]$.  
 There are many literature available related to decomposition of a continuous function into sum of two continuous functions, where each graph has pre-decided fractal dimensions. See \cite{hyde,humke,balka} for more details.
 \par  Motivated from above results, our aim is to study fractal dimensional results for the graph of product of two continuous functions.   
 In this article, we deal with this type of problems. Here we give  formulas for upper box dimension and lower box dimension of the graph of product of two continuous functions. We  give some results related to decomposition of continuous function into product of two continuous functions in terms of upper box dimension. We also establish an upper bound for upper box dimension of every element of a ring of polynomials in finite number of continuous functions over the field $\mathbb{R}$ and give some decomposition results for a continuous function.  
\par
   This paper is arranged as follows. In forthcoming Section \ref{2}, we give some basic definitions and required results for upcoming section. In Section \ref{3}, we characterise Theorem \ref{upf} and Theorem \ref{liu@} for product of two continuous functions and give some decomposition results for a continuous function. In last Section, we concluded our main results and discuss some open problems in this direction.      
	
	\section{Preliminaries}\label{2}
	\begin{defi}
		Let (X,d) be a metric space and $U$ be any nonempty subset of X. The diameter of $U$ is defined as $|U|=\sup\{d(x,y): x,y \in U\}$.
		Let $F$ be a subset of X. We say that $\{U_{i}\}$ is a $\delta$-cover of $F$ if $F\subset \bigcup\limits_{i=1}^{\infty}U_{i}$ with $0 <|U_{i}|\leq \delta$, for each ${i}$.
		\\ Let s be a non negative real number. For any $\delta>0$ we define
		$$H_{\delta}^{s}(F)=\inf\Big\{\sum_{i=1}^{\infty}|U_{i}|^{s} : \{U_{i}\} \text{ is a }\delta\text{-cover of }F\Big\}.$$
		As $\delta$ decreases, the class of permissible covers of $F$ reduces. Therefore, $H_{\delta}^{s}(F)$ increases, and so approaches to a limit, as $\delta\to 0$. We write
		$$H^{s}(F)=\lim_{\delta\to 0}H_{\delta}^{s}(F).$$
		We call $H^{s}(F)$ the $s$-dimensional Hausdorff measure of $F$.	
	\end{defi}
	\begin{defi}
		Hausdorff dimension of $F$ is the critical value of $s$ at which  $H^{s}(F)$ `jumps' from $\infty$ to $0$ and it is denoted by $\dim_{H}(F)$. Formally,
		$$\dim_{H}(F)=\inf\{s\geq 0 : H^{s}(F)=0\}=\sup\{s \geq 0 : H^{s}(F)=\infty\}$$
	\end{defi}
	\begin{defi}
		Let $F$ be any non-empty bounded subset of $X$ and let $N_{\delta}(F)$ be the smallest number of sets of diameter at most $\delta$ which can cover $F$. The lower and upper box-counting dimensions of $F$, respectively, are defined as
		$$\underline{\dim}_{B}F=\underline\lim_{\delta\to 0}\frac{\log{N_{\delta}(F)}}{-\log\delta},$$
		$$\overline{\dim}_{B}F=\overline\lim_{\delta \to 0}\frac{\log{N_{\delta}(F)}}{-\log\delta}.$$
		If these are equal we refer to the common value as the box-counting dimension or box dimension of $F$,
		$$\dim_{B}F=\lim_{\delta \to 0}\frac{\log{N_{\delta}(F)}}{-\log\delta}.$$
	\end{defi}
\begin{defi}
	Let $(X,d)$ be a metric space and $F\subset X$. For $\delta>0$ and $t\geq 0$, we define 
	$\mathcal{P}^t_{\delta}(F)=\sup\big\{\sum\limits_{i}|V_i|\big\}$, where $ \{V_i\}$ is a collection of disjoint balls of radii at most $\delta $ with centres in $F$.  As $\delta $ decreases, $\mathcal{P}^t_{\delta}(F)$ also decreases. Therefore, the limit
	   $$\mathcal{P}^t_{0}(F)=\lim_{\delta \to 0}\mathcal{P}^t_{\delta}(F)$$
	   exists. By the help of above limit, we define  
	   $$\mathcal{P}^t{(F)}= \inf\bigg\{\sum_{i}\mathcal{P}^t_{0}(F_i): F\subset \cup^{\infty}_{i=1}{F_i}\bigg\},$$ which is known as the $t$-dimensional packing measure. Therefore, packing dimension is defined as follows
	     $$\dim_{P}(F)=\inf\{t\geq 0 : \mathcal{P}^{t}(F)=0\}=\sup\{t\geq 0 : \mathcal{P}^{t}(F)=\infty\}$$
\end{defi}
\begin{defi}
 The Assouad dimension of a non-empty set $F\subseteq \mathbb{R}^2 $ is defined as
 \begin{align*}
 	 \dim_A(F) = \inf \bigg\{ \alpha: \exists~C >0~~\text{such that}~~\forall~ 0<r<R~~\text{and}~~ x\in F,~~ \\ N_r(B(x,R)\cap F)\leq C \bigg (\frac{R}{r}\bigg)^\alpha  \bigg\},
 \end{align*}
where $B(x,R)$ denotes the ball of radius $R$ with centre $x$ and $N_r(B(x,R)\cap F)$ is the  number of squares of $r $- mesh that intersect with $B(x,R)\cap F.$

\end{defi}

\begin{prop}\cite{Falconer}\label{p1}
	Let $f:[0,1]\to \mathbb{R}$ be a continuous function. Suppose that   $ \delta \in(0,1) $ and $m$ is the smallest natural number greater or equal to $\frac{1}{\delta }.$ Then, if $N_\delta(G(f))$ is the number of squares of $\delta $- mesh that intersect with graph f,
	$$\delta ^{-1}\sum_{i=0}^{m-1} R_{f}[i\delta,(i+1)\delta]\leq N_\delta(G(f)) \leq 2m+ \delta ^{-1}\sum_{i=0}^{m-1} R_{f} [i\delta ,(i+1)\delta], $$
	where  $R_f[i\delta,(i+1)\delta] = \sup\limits_{x,y\in [i\delta,(i+1)\delta]}|f(x)-f(y)|.$
	
\end{prop}

\begin{lemma}\cite{liu03} \label{LF}
	Let $f,g\in C[0,1]$. Then 
	$$\underline\dim_B{G{(f+g)}}\leq \max \{ \underline\dim_B{G(f)},\overline\dim_B{G(g)} \}.$$
\end{lemma}

\begin{prop}\cite{bayar}\label{DH}
Let $X$ be a compact subset of $\mathbb{R}$	with $\dim_{H}(X)>0$. Then every continuous functuon $f$ on $X$ can be decomposed as a sum of two continuous functions $g,h $ on $X$ such that 
 $$f=g+h~~\text{and}~~\dim_{H}G(g)=\dim_{H}G(h)=\dim_{H}(X) +1.$$
\end{prop}

\begin{theorem}\label{CB}
	Let $X$ be a compact subset of $\mathbb{R}$. Then for any $f\in C(X)$, we have
	$$\overline{\dim}_B(X)\leq \overline{\dim}_BG_f(X)\leq \overline{\dim}_B(X)+1,$$
	$$\underline{\dim}_B(X)\leq \underline{\dim}_BG_f(X)\leq \underline{\dim}_B(X)+1.$$
\end{theorem}
\begin{pff}
	One can easily prove this proposition by using monotonicity and Lipschitz invariance property of the upper and lower box dimension.
\end{pff}
\begin{prop}\cite{liu03}\label{EB}
	Let $X$ be a compact subset of $[0,1]$. Then for each continuous function $f$ on $X$, we have 
	$$\overline{\dim}_BG(F)\leq \overline{\dim}_B(X)+ 1,$$
	$$\underline{\dim}_BG(F)\leq \underline{\dim}_B(X)+ 1,$$
	where $F$ is the continuous linear extension of $f$ on $[0,1].$ 
\end{prop}
	\section{Main Results}\label{3}
	
 In the following lemma, we obtain a  formula for the upper box dimension of graph of product of two continuous functions.  
	\begin{lemma}\label{p2} Let $f,g\in C[0,1]$. Then 
		$$\overline{\dim}_{B}G{(f\cdot g)}\leq \max{\big\{\overline{\dim}_{B}G{(f)},\overline{\dim}_{B}G{(g)}\big \}}.$$	
	\end{lemma}
	\begin{pff} Let $f,g \in C[0,1],$
	\begin{align*}
		R_{f.g}[i\delta,(i+1)\delta]&=\sup_{x,y\in [i\delta,(i+1)\delta]}|(f.g)(x)-(f.g)(y)|\\
		&=\sup_{x,y\in [i\delta,(i+1)\delta]}|f(x)g(x)-f(x)g(y)+f(x)g(y)-f(y)g(y)|\\
		&=\sup_{x,y\in [i\delta,(i+1)\delta]}|f(x)(g(x)-g(y))+(f(x)-f(y))g(y)|\\
		&\leq \sup_{x,y\in [i\delta,(i+1)\delta]}|f(x)(g(x)-g(y))| +\sup_{x,y\in [i\delta,(i+1)\delta]}|(f(x)-f(y))g(y)|.	\end{align*}
Hence
\begin{equation}\label{e1}
	 R_{f\cdot g}[i\delta,(i+1)\delta]\leq M_1R_g[i\delta,(i+1)\delta]+M_2R_f[i\delta,(i+1)\delta],
\end{equation}	
where $M_1=\sup\limits_{x\in[0,1]}|f(x)| $ and $M_2=\sup\limits_{x\in[0,1]}|g(x)|.$
  Set $s=\max\{\overline{\dim}_{B}G_{f},\overline{\dim}_{B}G_{g}\big \}.$ Let $\epsilon>0,$ then by the definition of the upper box dimension there exists a $\delta_0\in (0,1)$ such that for all $\delta\leq \delta_0$, we have
  \begin{align*}
  	N_\delta(G(f))&\leq \delta^{-s-\epsilon},\\
  		N_\delta(G(g))&\leq \delta^{-s-\epsilon}.
  \end{align*}
   Using the above inequalities and Proposition \ref{p1} , we get
\begin{equation}\label{e2}
		\delta ^{-1}\sum_{i=0}^{m-1} R_{f}[i\delta,(i+1)\delta]\leq \delta^{-s-\epsilon},
\end{equation}
and	
\begin{equation}\label{e3}
	\delta ^{-1}\sum_{i=0}^{m-1} R_{g}[i\delta,(i+1)\delta]\leq \delta^{-s-\epsilon}.
\end{equation}
Therefore by inequalities \eqref{e1}, \eqref{e2} and \eqref{e3}, we have 
\begin{equation}
\begin{aligned}
	\delta^{-1}	\sum_{i=0}^{m-1} R_{f.g}[i\delta,(i+1)\delta]&\leq M_1	\delta^{-1}\sum_{i=0}^{m-1} R_g[i\delta,(i+1)\delta]+M_2	\delta^{-1}\sum_{i=0}^{m-1} R_f[i\delta,(i+1)\delta]\\
	&\leq M_1\delta^{-s-\epsilon}+M_2\delta^{-s-\epsilon}\\
	&=C\delta^{-s-\epsilon},
\end{aligned}	
\end{equation}
	
where $C=M_1+M_2.$ 	
Hence by above and Proposition \ref{p1} , we have	
\begin{equation}\label{e4}
\begin{aligned}
	N_\delta(G(f\cdot g))& \leq 2m+ C\delta^{-s-\epsilon}\\
	&\leq 2(1+\delta^{-1})+ C\delta^{-s-\epsilon}\\
	&\leq (4+C)\delta^{-s-\epsilon}\\
	&=C_1\delta^{-s-\epsilon}, 
\end{aligned}	
\end{equation}
where $C_1=C+4$. From inequality \ref{e4}, we conclude that
$$\overline{\dim}_{B}G(f\cdot g)=\varlimsup_{\delta \rightarrow 0} \frac{\log N_{\delta}(G(f\cdot g))}{- \log \delta}\leq s+\epsilon.$$
Since this is true for all $\epsilon>0$,
$$\overline{\dim}_{B}G{(f\cdot g)}\leq \max{\big\{\overline{\dim}_{B}G{(f)},\overline{\dim}_{B}G{(g)}\big \}}.$$
This completes the proof.
	\end{pff}
  
  The following example show that the inequality in the previous lemma may be strict.
  \begin{example}
  	let $g:[0,1]\to [0,1]$ be component of Peano space filling curve \cite{kono} and $f:[0,1]\to \mathbb{R} $ such that $f(x)=0$ for all $x\in [0,1].$
  	 Then it is well- known that $\overline{\dim}_BG(g)=1.5$, $\overline{\dim}_BG(f)=1$ and  $\overline{\dim}_BG(f\cdot g)=1.$

  \end{example}

Next lemma provides a dimensional relation between a continuous function and its reciprocal function.
\begin{prop}\label{l3.3}
If $f\in C[0,1]$ such that $f(x)\ne 0$ for all $x\in[0,1]$. Then 
\begin{align*}
	\overline{\dim}_BG(f)=&\overline{\dim}_BG\bigg(\frac{1}{f}\bigg),~~  \underline{\dim}_BG(f)=\underline{\dim}_BG\bigg(\frac{1}{f}\bigg),~~ {\dim}_PG(f)={\dim}_PG\bigg(\frac{1}{f}\bigg), \\&{\dim}_HG(f)={\dim}_HG\bigg(\frac{1}{f}\bigg)~~\text{and}~~  {\dim}_AG(f)={\dim}_AG\bigg(\frac{1}{f}\bigg).
\end{align*}

 	\end{prop}	
 \begin{pff} Suppose that $M_1=\inf\limits_{x\in [0,1]}|f(x)|$ and $M_2=\sup\limits_{x\in[0,1]}|f(x)|$.
 	Let us define a mapping $\Psi : G(f) \to   G\bigg(\frac{1}{f}\bigg)$ by $$\Psi(x,f(x))=\bigg(x,\bigg(\frac{1}{f}\bigg)(x)\bigg).$$ Our aim is to show that $\Psi$ is a bi-Lipschitz map. Using the simple properties of norm, it follows that
 	\begin{align*}
 		\|\Psi(x,f(x))-\Psi(y,f(y))\|^2&=\bigg\|\bigg(x,\bigg(\frac{1}{f}\bigg)(x)\bigg)-\bigg(y,\bigg(\frac{1}{f}\bigg)(y)\bigg) \bigg\|^2 \\
 		&=|x-y|^2 +\bigg|\bigg(\frac{1}{f}\bigg)(x)\bigg)-\bigg(\frac{1}{f}\bigg)(y)\bigg)\bigg|^2\\
 		&=|x-y|^2+\bigg|\frac{f(x)-f(y)}{f(x)f(y)}\bigg|^2 \\
 		&\leq |x-y|^2+\bigg(\frac{1}{{M_1}^4}\bigg)|{f(x)-f(y)}\big|^2\\
 		&\leq \bigg(1+\bigg(\frac{1}{{M_1}^4}\bigg)\bigg) \bigg(|x-y|^2+ |{f(x)-f(y)}\big|^2\bigg)\\
 		&= C_1 \|(x,f(x))-(y,f(y))\|^2,
 		\end{align*}
  where $C_1=\bigg(1+\bigg(\frac{1}{{M_1}^4}\bigg)\bigg).$ Also,
  \begin{align*}
  	\|(x,f(x))-(y,f(y))\|^2&= |x-y|^2+ |{f(x)-f(y)}\big|^2\\
  	&=|x-y|^2+\bigg|\frac{f(x)-f(y)}{f(x)f(y)}\bigg|^2|f(x)f(y)|^2\\
  	 &\leq|x-y|^2 +\bigg|\bigg(\frac{1}{f}\bigg)(x)\bigg)-\bigg(\frac{1}{f}\bigg)(y)\bigg)\bigg|^2{M_2}^4\\
  	 & \leq (1+{M_2}^4)\bigg(|x-y|^2 +\bigg|\bigg(\frac{1}{f}\bigg)(x)\bigg)-\bigg(\frac{1}{f}\bigg)(y)\bigg)\bigg|^2\bigg)\\
  	 &= C_2 \|\Psi(x,f(x))-\Psi(y,f(y))\|^2,
  	 \end{align*}
   where $C_2=(1+{M_2}^4)$. By the above two inequalities, we get 
   $$ {C_2}^{-1} \|(x,f(x))-(y,f(y))\|^2\leq	\|\Psi(x,f(x))-\Psi(y,f(y))\|^2\leq C_1 \|(x,f(x))-(y,f(y))\|^2.$$
   Therefore, $\Psi$ is a bi-Lipschitz map. By using the bi-Lipschitz invariance  property of Hausdorff dimension, upper box dimension, lower box dimension, packing dimension and Assouad dimension, we get our required result.
    \end{pff}

 In the upcoming lemma, we describe a relation between a continuous function and its square function for various dimensions.
\begin{prop}\label{l3.4}
	Let $f\in C[0,1]$ be such that $f(x)\ne 0$ for all $x\in [0,1]$. Then 
	\begin{align*}
		\overline{\dim}_BG(f^2)=&\overline{\dim}_BG(f),~~\underline{\dim}_BG(f^2)=\underline{\dim}_BG(f),~~{\dim}_PG(f^2)={\dim}_PG(f),\\&{\dim}_HG(f^2)={\dim}_HG(f)~~\text{and}~~{\dim}_AG(f^2)={\dim}_AG(f).
	\end{align*}
	
\end{prop}	
\begin{pff}  Without loss of generality, we can assume that $f(x)>0$ for all $x \in [0,1]$ and  suppose $M_1=\inf\limits_{x\in [0,1]}f(x)$ and $M_2=\sup\limits_{x\in[0,1]}f(x)$.
	Let us define a mapping $\Phi : G(f) \to   G(f^2)$ as follows
	$$\Phi\big(x,f(x)\big)=\big(x,f^2(x)\big).$$
	We claim that $\Phi$ is a bi-Lipschitz map. Using the simple properties of norm, it follows that
	 \begin{align*}
	 	\|\Phi(x,f(x))-\Phi(y,f(y))\|^2&= \|(x,f^2(x))-(y,f^2(y))\|^2\\
	 	&=|x-y|^2 +|f^2(x)-f^2(y)|^2\\
	 	&= |x-y|^2+ |f(x)-f(y)|^2 |f(x)+f(y)|^2\\
	 	&\leq |x-y|^2+ 4 {M_2}^2|f(x)-f(y)|^2 \\
	 	&\leq (1+4 {M_2}^2) \big(|x-y|^2+|f(x)-f(y)|^2\big)\\
	 	&= C_3 \|(x,f(x))-(y,f(y))\|^2,
	 \end{align*}
where $C_3=(1+4 {M_2}^2).$ Also,
\begin{align*}
	\|(x,f(x))-(y,f(y))\|^2&=|x-y|^2+|f(x)-f(y)|^2\\
	&= |x-y|^2+ \frac{|f^2(x)-f^2(y)|^2}{|f(x)+f(y)|^2}\\
	&\leq |x-y|^2+\bigg(\frac{1}{4 {M_1}^2}\bigg) |f^2(x)-f^2(y)|^2\\
	&\leq \bigg(1+\bigg(\frac{1}{4 {M_1}^2}\bigg)\bigg)\bigg(|x-y|^2+|f^2(x)-f^2(y)|^2\bigg)\\
	&= C_4 \|\Phi(x,f(x))-\Phi(y,f(y))\|^2,
	\end{align*}
where $C_4=\bigg(1+\bigg(\frac{1}{4 {M_1}^2}\bigg)\bigg).$ Using the above two inequalities, we get	
$$ {C_4 }^{-1}	\|(x,f(x))-(y,f(y))\|^2\leq \|\Phi(x,f(x))-\Phi(y,f(y))\|^2\leq C_3 \|(x,f(x))-(y,f(y))\|^2.$$
Therefore, $\Phi$ is a bi-Lipschitz map. We get our required result, by using the bi-Lipschitz invariance property of Hausdorff dimension, upper box dimension, lower box dimension, packing dimension and Assouad dimension. 
	
\end{pff}

Next, we give a generalized version of the above result.
\begin{prop}
	Let $f\in C[0,1]$ be such that $f(x)\ne 0$ for all $x\in [0,1]$. Then for any $ n\in \mathbb{N} $, we have
	\begin{align*}
		\overline{\dim}_BG(f^n)=&\overline{\dim}_BG(f),~~\underline{\dim}_BG(f^n)=\underline{\dim}_BG(f),~~{\dim}_PG(f^n)={\dim}_PG(f),\\&{\dim}_HG(f^n)={\dim}_HG(f)~~\text{and }~~{\dim}_AG(f^n)={\dim}_AG(f).
	\end{align*}
\end{prop}

    \begin{pff}
      One can easily prove this proposition by following the proof of Proposition \ref{l3.4}.
\end{pff}

In the following proposition, we establish an equality relation for upper box dimension of graph of product of two continuous functions. 
\begin{prop}\label{l3.5}
	Suppose that $f,g\in C[0,1]$ are such that $f(x)\ne 0, g(x)\ne 0$ for all $x\in [0,1]$ and 
	$$\overline{\dim}_B G(f)\ne \overline{\dim}_BG(g).$$
	Then  
	$$\overline{\dim}_{B}G{(f\cdot g)}= \max{\big\{\overline{\dim}_{B}G{(f)},\overline{\dim}_{B}G{(g)}\big \}}.$$	
\end{prop}	
\begin{pff}
	Without loss of generality, we can assume that
	\begin{equation}\label{3.51}
		\overline{\dim}_B G(f)<\overline{\dim}_BG(g).
	\end{equation}
	Let $h=f\cdot g$ and assume that	 $$\overline{\dim}_BG(h)\ne\overline{\dim}_BG(g).$$
	 By using Lemma \ref{p2} and by the above assumption, we have
	\begin{equation}\label{3.52}
		\overline{\dim}_BG(h)<\overline{\dim}_BG(g).
	\end{equation}
	By inequalities \eqref{3.51} and \eqref{3.52}, we get 
	$$\overline{\dim}_BG\bigg(h\cdot \frac{1}{f}\bigg)=\overline{\dim}_BG(g)>\max{\big\{\overline{\dim}_{B}G{(h)},\overline{\dim}_{B}G{(f)}\big \}}.$$
	Hence by Proposition \ref{l3.3} and the above inequality, we get
	$$\overline{\dim}_BG\bigg(h\cdot \frac{1}{f}\bigg)>\max{\bigg\{\overline{\dim}_{B}G{(h)},\overline{\dim}_{B}G\bigg(\frac{1}{f}\bigg)\bigg \}},$$
	which is a contradiction to Lemma \ref{p2}. Therefore,
	$$\overline{\dim}_{B}G{(f\cdot g)}= \max{\big\{\overline{\dim}_{B}G{(f)},\overline{\dim}_{B}G{(g)}\big \}}.$$ 
	This completes the proof.
\end{pff}

In the above proposition the condition $\overline{\dim}_B G(f)\ne \overline{\dim}_BG(g)$ can not be dropped. For the support of this we give the following example.
\begin{example}
	For any $s\in (1,2)$, we can define a Weierstrass type function $f:[0,1]\to \mathbb{R}$ such that $f(x)\ne 0$ for all $x\in [0,1]$ with the property $\overline{\dim}_BG(f)=s.$ Take $g:[0,1]\to \mathbb{R}$ such that $g(x)=\frac{1}{f(x)}$ for all $x\in [0,1].$ Therfore by Proposition \ref{l3.3}, $\overline{\dim}_BG(g)=s.$ 
\end{example}

In the upcoming proposition, we obtain a formula for lower box dimension of the graph of product of two continuous functions.
\begin{prop}\label{luff}
	Let $f,g\in C[0,1].$ Then
	$$\underline{\dim}_{B}G{(f\cdot g)}\leq \max{\big\{\underline{\dim}_{B}G{(f)},\overline{\dim}_{B}G{(g)}\big \}}.$$ 
\end{prop}
\begin{pff} Let $M_1=\sup\limits_{x\in[0,1]}|f(x)| $ and $M_2=\sup\limits_{x\in[0,1]}|g(x)|. $ Then, as in inequality \ref{p2},
	\begin{equation}\label{e3.16}
		R_{f.g}[i\delta,(i+1)\delta]\leq M_1R_g[i\delta,(i+1)\delta]+M_2R_f[i\delta,(i+1)\delta],	
	\end{equation}
	Let $\delta = 5^{-n}$ and $\alpha = \max{\big\{\underline{\dim}_{B}G{(f)},\overline{\dim}_{B}G{(g)}\big \}}.$
	\\ By the definition of $ \underline{\dim}_{B}G{(f)}$, there exists a subsequence $\{n_i\}$ of natural numbers such that 
	$$\lim_{i\to \infty }\frac{\log N_{5^{-n_i}}G(f)}{-\log 5^{-n_i}}=\underline{\dim}_{B}G{(f)}.$$
	Let $\epsilon >0,$ then there exists $i_0$ such that 
	\begin{equation}\label{e3.17}
		N_{5^{-n_i}}G(f)\leq 5^{n_i(\alpha +\epsilon) }, 
	\end{equation}
	for all $i\geq i_0.$  By the definition of $\overline{\dim}_{B}G{(g)}$, there exists  $N\in \mathbb{N}$ such that 
	\begin{equation}\label{e3.18}
		N_{5^{-n}}G(g)\leq 5^{n(\alpha +\epsilon) },
	\end{equation}  
	for all $n\geq N$.
	Using Proposition \ref{p1} and inequality \eqref{e3.16}, we get
	$$N_{5^{-n_i}}G(f\cdot g)\leq 2m+ N_{5^{-n_i}}G(f)+N_{5^{-n_i}}G(g)~~~\text{for all}~~i\geq1.$$
	Therefore by inequalities \eqref{e3.17} and \eqref{e3.18}, we have  
	\begin{align*}
		N_{5^{-n_i}}G(f\cdot g)&\leq 2m+ 5^{n_i(\alpha +\epsilon) }+5^{n_i(\alpha +\epsilon) }\\
		&\leq 6 \cdot 5^{n_i(\alpha +\epsilon) },
	\end{align*}
	for all $n_i\geq\max \{n_{i_0},N\}$. Then
	\begin{align*}
		\underline{\dim}_BG(f\cdot g)&\leq \lim_{i\to \infty }\frac{\log N_{5^{-n_i}}G(f.g)}{-\log 5^{-n_i}} \\
		&\leq \alpha +\epsilon.
	\end{align*}
	Since $\epsilon$ is arbitrary, we have
	$$\underline{\dim}_{B}G{(f\cdot g)}\leq \max{\big\{\underline{\dim}_{B}G{(f)},\overline{\dim}_{B}G{(g)}\big \}}.$$ 
	This completes the proof.
\end{pff}

In the following theorem, we establish a general upper bound for upper box dimension of a graph of any  polynomial in two continuous functions over the field $\mathbb{R}.$

 \begin{theorem}\label{ply}
 	Let $f,g\in C[0,1]$ and let $\mathcal{R}$ denote the ring of polynomials $P(f,g)$ in  $f,g$ over the field $\mathbb{R}.$ Then 
 	$$\overline{\dim}_BG(P(f,g))\leq \max{\big\{\overline{\dim}_{B}G{(f)},\overline{\dim}_{B}G{(g)}\big \}},$$
 	
 	for any $P(f,g)\in \mathcal{R}$.
 	 	
 \end{theorem}
 \begin{pff} Let $ f,g\in C[0,1]$. Theorem \ref{upf} yields,
 	\begin{equation}\label{3.11}
 		\overline\dim_B{G{(f+g)}}\leq \{ \overline\dim_B{G(f)},\overline\dim_B{G(g)} \}.
 	\end{equation}

 By using Lemma \ref{p2}, we get
 \begin{equation}\label{3.13}
 	\overline{\dim}_BG(f^n)\leq \overline{\dim}_BG(f),
 \end{equation}
for each $n\in \mathbb{N}.$
Combining Lemma \ref{p2}, inequalities \eqref{3.11} and \eqref{3.13}, we get the required result.  This completes the proof. 
 \end{pff}

The next result generalizes the above theorem for any finte number of continuous functions.
\begin{prop}
	If $f_i\in C[0,1]$ for $i\in \{1,2,\cdots,m\}$ and $\mathcal{R}^{'}$ denotes the ring of polynomials in $f_i\in C[0,1]$ for $i\in \{1,2,\cdots,m\}$ over the field $\mathbb{R}.$ Then 
	$$\overline{\dim}_BG(r^{'})\leq \max{\big\{\overline{\dim}_{B}G{(f_1)},\overline{\dim}_{B}G{(f_2)},\cdots ,\overline{\dim}_{B}G{(f_m)} \big \}},$$
	 for any $r^{'}\in \mathcal{R}^{'}$.
\end{prop}
\begin{pff}
	By using the similar arguments in the proof of Theorem \ref{ply}.
\end{pff}

In the next proposition, we obtain a general upper bound for upper box dimension of the graph of  any rational function in two continuous functions over the field $\mathbb{R}$.
\begin{prop}\label{pr2}
	Let $f,g\in C[0,1]$ and let $\mathcal{Q}$ denote the ring of rational functions  $R(f,g)$ over the field $\mathbb{R},$ where $$R(f,g)=\frac{P(f,g)}{Q(f,g)},$$ 
	where $P(f,g)$ and $Q(f,g)$ are polynomials in $f,g$ over the field $\mathbb{R}$ and $Q(f,g)(x)\ne 0$ for any $x\in [0,1].$ Then
	$$\overline{\dim}_BG(R(f,g))\leq \max{\big\{\overline{\dim}_{B}G{(f)},\overline{\dim}_{B}G{(g)}\big \}}$$
	
	for any $ R(f,g) \in \mathcal{Q}$.
	
\end{prop} 
 \begin{pff}
By using Proposition \ref{l3.3}, we get
\begin{equation}\label{eqi}
\overline\dim_{B}G(Q(f,g))=\overline{\dim}_BG\bigg(\frac{1}{Q(f,g)}\bigg).	
\end{equation}

In the light of Theorem \ref{ply}, Lemma \ref{p2} and equation \eqref{eqi}, we obtain our required result.
 \end{pff}

The next result generalizes the previous proposition for any finte number of continuous functions.

\begin{prop}
	Let $f_i\in C[0,1]$ for $i\in \{1,2,\cdots,m\}$ and let $\mathcal{Q}^{'}$ denote the ring of rational functions  $R(f_1,f_2,\cdots,f_m)$ over the field  $\mathbb{R},$ where $$R(f_1,f_2,\cdots,f_m)=\frac{P(f_1,f_2,\cdots,f_m)}{Q(f_1,f_2,\cdots,f_m)},$$ 
	where $P(f_1,f_2,\cdots,f_m)$ and  $Q(f_1,f_2,\cdots,f_m)$ are polynomials in $f_i$ for $i\in \{1,2,\cdots,m\}$ over the field $\mathbb{R}$ and $Q(f_1,f_2,\cdots,f_m)(x)\ne 0$ for any $x\in [0,1].$ Then
	$$\overline{\dim}_BG(q')\leq \max\limits_{i\in \{1,2,\cdots,m\} }{\big\{\overline{\dim}_{B}G{(f_i)}\big \}}$$
	
	for any $q'\in \mathcal{Q}^{'}$.
\end{prop}
\begin{pff}
	Idea of proof of Proposition \ref{pr2} yields our result. 
\end{pff}

In the following theorems, we decompose a continuous function into product of two continuous with some specific dimensional properties.  
 	
\begin{theorem}\label{MD}
	Let $\beta \in [1,2]$ and let  $f\in C[0,1]$ be a given function such that $f(x)\ne 0$ for all $x\in [0,1]$. Then there exist two functions $g,h\in C[0,1]$ such that 
	$$ f=g\cdot  h~~~\text{and}~~~\overline{\dim}_BG(g)=\overline{\dim}_BG(h)=\beta $$ if and only if  $~~\overline{\dim}_BG(f)\leq \beta.$ 
	
\end{theorem}
\begin{pff} The `only if' part follows directly from Lemma \ref{p2}. Hence we only need to prove the `if' part.\par 
	Let  $\overline{\dim}_BG(f)= \beta$. If we take $g=\frac{1}{f}$ and $h=f^2$. Then, by Proposition \ref{l3.3} and Proposition \ref{l3.4}, we get
	$$ f=g\cdot  h~~~~\text{and}~~~~\overline{\dim}_BG(g)=\overline{\dim}_BG(h)=\beta. $$ 
	Let   $\overline{\dim}_BG(f)< \beta$.
	Choose a non-empty perfect set $X\subset [0,1]$ such that $\dim_{H}(X)=\overline{\dim}_{B}(X)=\beta -1.$ By Proposition \ref{DH}, there exists a function $f_1\in C(X)$ such that 
	$$\dim_{H}G_{f_1}(X)=\beta.$$
	 By using Theorem \ref{CB}, we have
	$$\dim_{H}G_{f_1}(X)=\overline{\dim}_BG_{f_1}(X)=\beta.$$ 
	 Let $k$ be the continuous linear extension of $f_1$ on $[0,1].$ By the countable stability property of Hausdorff dimension and Proposition \ref{EB}, we get
	  $$\dim_{H}G(k)=\overline{\dim}_BG(k)=\beta.$$
	  By adding a suitable constant, we can take $k$ such that $k(x)\ne 0$ for all $x\in [0,1].$

	Let $g=k\cdot f$. So by Proposition \ref{l3.5}, 
	$$\overline{\dim}_BG(g)=\max{\big\{\overline{\dim}_{B}G{(k)},\overline{\dim}_{B}G{(f)}\big \}}=\beta.$$
	Let $h=\frac{1}{k}$. Therefore by Proposition \ref{l3.3}, we have
	$$\overline{\dim}_BG(h)=\overline{\dim}_BG(k)=\beta.$$
	Hence $$f=g\cdot h~~~\text{and}~~~ \overline{\dim}_BG(g)=\overline{\dim}_BG(h)=\beta. $$
	This completes the proof.
	\end{pff}
\begin{theorem}
	Let $f\in C[0,1]$ be such that $f(x)\ne 0$ for all $x\in [0,1]$ and $\overline{\dim}_BG(f)=\alpha .$ Then for any $\beta \in [1,\alpha)$,
	there exist two continuous functions $g$ and $h$ in $C[0,1]$ such that
	$$f=g\cdot h,~~\overline{\dim}_B G(g)=\alpha~~~\text{and}~~~\overline{\dim}_BG(h)=\beta.$$
\end{theorem}
\begin{pff}
By using the same technique as in Theorem \ref{MD}, we can choose a function $k\in C[0,1]$ such that $k(x)\ne 0$ for all $x\in C[0,1]$ and $\overline{\dim}_BG(k)=\beta.$
Let $g=k\cdot f$. Then by using  Proposition \ref{l3.5}, we get $$\overline{\dim}_BG(g)=\alpha.$$
Let $h=\frac{1}{k}$. Now, by Proposition \ref{l3.3}, we get $$\overline{\dim}_BG(h)=\beta.$$
So, we conclude that
	$$f=g\cdot h,~~\overline{\dim}_B G(g)=\alpha~~~\text{and}~~~\overline{\dim}_BG(h)=\beta.$$
	This completes the proof.
	\end{pff}
\section{Conclusion and some open problems}	
\par 
Firstly, we gave an upper bound for the upper box dimension of the graph of product of two continuous functions in terms of their upper box dimension, which is an analogue of Theorem \ref{upf}. We also establish an upper bound for lower box dimension of the graph of product of two continuous functions. Next, we prove that we can decompose a continuous function as a product of two continuous functions where upper box dimension of the graph of these continuous functions is same as any number which is greater than or equal to the upper box dimension of graph of the given continuous function.
\par
In this paper, we proved the decomposition results for the upper box dimension in terms of product, which is a analogue of Theorem \ref{liu@}, but the following analogous decomposition results regarding Hausdorff dimension \cite[Theorem 1.2]{liu01}, packing dimension \cite[Theorem 1.7]{liu02} and lower box dimension \cite[Theorem 1.4]{liu03} are still open.
  \begin{quest} 
  	  For any $\beta \in [1,2]$ and a given function $f\in C[0,1]$, does there exist two continuous functions $g,h\in C[0,1]$ such that
  	  $$f=g\cdot h~~\text{and}~~~ \dim_{H}G(g)=\dim_{H}G(h)=\beta? $$
\end{quest}

\begin{quest}
	Let $\beta\in [1,2]$ and $f\in C[0,1]$. If $\dim_PG(f)\leq \beta$, does there exist functions $g,h\in C[0,1]$ such that
	$$f=g\cdot h~~\text{and }~~\dim_{P}G(g)=\dim_{P}G(h)=\beta?$$
\end{quest}
\begin{quest}
 For any $\beta \in [1,2]$ and a given function $f\in C[0,1]$, does there exist two continuous functions $g,h\in C[0,1]$ such that
$$f=g\cdot h~~\text{and}~~~ \underline{\dim}_{B}G(g)=\underline{\dim}_{B}G(h)=\beta?$$
\end{quest}
\bibliographystyle{amsplain}

\begin{thebibliography}{10}
		\bibitem{balka} R. Balka, Dimensions of graphs of prevalent continuous maps, Journal of Fractal Geometry, 3(4) (2016) 407-428.
	\bibitem{bayar}	F. Bayart, Y. Heurteaux, On the Hausdorff dimension of graphs of prevalent continuous functions on compact sets, in: Further developments in fractals and related fields, Birkhauser, Boston, 2013, pp. 25-34.
	\bibitem{Falconer}  K. J. Falconer, Fractal Geometry: Mathematical Foundations and Applications, John Wiley Sons Inc., New York, 1999.
		\bibitem{falconer2011horizon} K. J. Falconer, J. M. Fraser, The horizon problem for prevalent surfaces, Mathematical Proceedings of the Cambridge Philosophical Society, 151(2) (2011) 355–372.
		\bibitem{fraser} J. M. Fraser, Assouad dimension and fractal geometry, Cambridge University Press. 2020.
			\bibitem{humke} PD. Humke, G. Petruska, The packing dimension of a typical continuous function is 2, Real Analysis Exchange, 14(2) (1988) 345-358.
	\bibitem{hyde} J. Hyde, V. Laschos, L. Olsen, I. Petrykiewicz, A. Shaw, On the box dimensions of graphs of typical continuous functions, Journal of Mathematical Analysis and Applications, 391(2) (2012) 567-581.
		\bibitem{kono} N. K{\^o}no, On self-affine functions, Japan Journal of Applied Mathematics, 3(2) (1986) 259-269.
\bibitem{liu01}	J. Liu, J. Wu, A remark on decomposition of continuous functions, Journal of Mathematical Analysis and Applications, 401 (2013) 404-406.
\bibitem{liu02} J. Liu, B. Tan, J. Wu, Graphs of continuous functions and packing dimension, Journal of Mathematical Analysis and Applications,
435(2) (2016) 1099-1106.
\bibitem{liu03} J. Liu, D. Liu, On the Decomposition of Continuous Functions and Dimensions, Fractals, 28(01) (2020) 2050007(6 pages)

\bibitem{mauldin1}	R. D. Mauldin, S. C. Williams, On the Hausddorff dimension of some graphs, Trans. Amer. Math. Soc. 298 (1986) 789-803.

	\bibitem{win} P. Wingren, Dimensions of graphs of functions and lacunary decompositions of spline approximations, Real Anal. Exchange (2000) 17-26.






\end{thebibliography}

\end{document}